\documentclass[12pt]{article}
\usepackage[latin9]{inputenc}
\usepackage{amsmath}
\usepackage{amssymb}
\usepackage{esint}

\makeatletter
\usepackage{mathtools}
\newtheorem{thm}{Theorem} [section]

 \newtheorem{rem}[thm]{Remark}     

\newcommand{\coloneqq}{:=}

\renewcommand{\eqref}[1]{(\ref{#1})}
\renewcommand{\overset}[2]{\stackrel{#1}{#2}}

\DeclareMathAccent{\Circ}{\mathalpha}{operators}{"17}
\newcommand{\interior}[1]{\Circ{#1}}

\renewcommand{\Re}{\operatorname{\mathfrak{Re}}}

\renewcommand{\tilde}{\widetilde}
\renewcommand*{\epsilon}{\varepsilon}
\renewcommand*{\rho}{\varrho}
\makeatother

\begin{document}

\title{On Boundary Damped Inhomogeneous Timoshenko Beams and Related Problems}

\author{R.~Picard%
\thanks{Institut f{ü}r Analysis,FR Mathematik, TU Dresden, D-01062 Dresden,
Germany%
} \, \& B.A.~Watson%
\thanks{School of Mathematics, University of the Witwatersrand, Johannesburg,
South Africa.%
}~ %
\thanks{Funded in part by NRF grant IF2011032400120 and the Centre for Applicable
Analysis and Number Theory%
}}
\maketitle
\begin{abstract}
We consider the model equations for the Timoshenko beam as a first
order system in the framework of evolutionary equations as developed
in \cite{Pi2009-1}. The focus is on boundary damping, which is implemented
as a dynamic boundary condition. A change of material laws allows
to include a large class of cases of boundary damping. By choosing
a particular material law, it is shown that the first order approach
to Sturm-Liouville problems with boundary damping is also covered. 
\end{abstract}

\section{Introduction\label{sec-intro}}

The homogeneous Timoshenko beam model is given by a second order system
of the form 
\begin{equation}
\begin{array}{rl}
\nu_{1}\partial_{0}^{2}\phi-\partial_{1}\kappa_{1}^{-1}\partial_{1}\phi+\kappa_{2}^{-1}\left(\partial_{1}u+\phi\right) & =f_{1}.\\
\nu_{2}\partial_{0}^{2}u-\partial_{1}\kappa_{2}^{-1}\left(\partial_{1}u+\phi\right)+d\partial_{0}u & =f_{2},
\end{array}\label{eq:Timoshenko-1}
\end{equation}
Here $u$ denotes the radial displacement of the beam and $\phi$
its angular displacement. The real time independent parameters $\nu_{1},\kappa_{1},\nu_{2},\kappa_{2},d$,
describe physical and geometrical properties of the beam, which, however,
we will not elaborate upon further, since for the analysis their interpretation
is not important. It may be interesting to note that these equations
are the 1-dimensional version of the Reissner-Mindlin plate model,
compare \cite{Picard201354}. To reformulate \eqref{eq:Timoshenko-1}
as a first order system we introduce the velocities $\eta\coloneqq\partial_{0}\phi$
and $s\coloneqq\partial_{0}u$ as new unknowns. Moreover, we let 
\begin{equation}
V_{2}\coloneqq\kappa_{2}^{-1}\left(\partial_{1}u+\phi\right)\label{eq:mod-Hooke}
\end{equation}
and 
\begin{equation}
V_{1}\coloneqq\kappa_{1}^{-1}\partial_{1}\phi.\label{eq:rot-Hooke}
\end{equation}
Now \eqref{eq:Timoshenko-1} can be written as 
\begin{eqnarray*}
\nu_{1}\partial_{0}\eta-\partial_{1}V_{1}+V_{2} & = & f_{1}\\
\nu_{2}\partial_{0}s-\partial_{1}V_{2}+ds & = & f_{2}
\end{eqnarray*}
together with the spatial differentiated equations \eqref{eq:mod-Hooke}
and \eqref{eq:rot-Hooke} resulting in the temporal equations 
\begin{eqnarray*}
\partial_{0}\kappa_{1}V_{1}-\partial_{1}\eta & = & 0,\\
\partial_{0}\kappa_{2}V_{2}-\partial_{1}s-\eta & = & 0.
\end{eqnarray*}
Written as a system we have 
\[
\left(\partial_{0}\mathcal{M}_{0}+\mathcal{M}_{1}+\left(\begin{array}{cccc}
0 & -\partial_{1} & 0 & 0\\
-\partial_{1} & 0 & 0 & 0\\
0 & 0 & 0 & -\partial_{1}\\
0 & 0 & -\partial_{1} & 0
\end{array}\right)\right)\left(\begin{array}{c}
V_{1}\\
\eta\\
s\\
V_{2}
\end{array}\right)=\left(\begin{array}{c}
0\\
f_{1}\\
f_{2}\\
0
\end{array}\right),
\]
where 
\[
\mathcal{M}_{0}\coloneqq\left(\begin{array}{cccc}
\kappa_{1} & \qquad0 & \qquad0 & \qquad0\\
0 & \qquad\nu_{1} & \qquad0 & \qquad0\\
0 & \qquad0 & \qquad\nu_{2} & \qquad0\\
0 & \qquad0 & \qquad0 & \qquad\kappa_{2}
\end{array}\right),\;\mathcal{M}_{1}\coloneqq\left(\begin{array}{cccc}
0 & \qquad0 & \qquad0 & \qquad0\\
0 & \qquad0 & \qquad0 & \qquad1\\
0 & \qquad0 & \qquad d & \qquad0\\
0 & \qquad-1 & \qquad0 & \qquad0
\end{array}\right)
\]
Here $\nu_{1},\kappa_{1},\nu_{2},\kappa_{2}$ are positive real time
independent parameters. It is essentially the latter form in which
we will approach the issue of boundary damping, since it will allow
us to utilize the framework of \cite{Pi2009-1}. The structural features
of the first order system allow us also to make -- almost effortlessly
-- the transition to more general media. %%%%%%%%%%%%%%%%%%%%%%%%%%%%%%%%%%%%%

\section{Functional Analytic Framework\label{sec-FA}}

Key to the approach presented here is the closure of the differentiation
operator on $C_{1}\left(\mathbb{R},H\right)$-functions with compact
support, i.e. acting on $\interior C_{1}\left(\mathbb{R},H\right)$.
It may assist the reader to consult \cite{Pi2009-1,PIC_2010:1889,PDE_DeGruyter,PIC-WAT}
for extra details on the approach used here. In particular we consider
differentiation as an operator in $H_{\rho,0}\left(\mathbb{R},H\right)$,
$\rho\in\left(0,\infty\right)$, a weighted $L^{2}$-type space with
inner product 
\[
\left\langle \varphi|\psi\right\rangle _{\rho,0,0}\coloneqq\int_{\mathbb{R}}\left\langle \varphi\left(t\right)|\psi\left(t\right)\right\rangle _{0}\:\exp\left(-2\rho t\right)\: dt,
\]
where $\left\langle \:\cdot\:|\:\cdot\:\right\rangle _{0}$ denotes
the inner product of $H$, which, as a matter of convention, we assume
to be linear in the second factor. The resulting operator 
\[
\partial_{0}:D\left(\partial_{0}\right)\subseteq H_{\rho,0}\left(\mathbb{R},H\right)\to H_{\rho,0}\left(\mathbb{R},H\right)
\]
turns out, \cite{PIC_2010:1889,PDE_DeGruyter}, to be normal. In particular,
$\partial_{0}^{*}=-\partial_{0}+2\rho$ and so 
\begin{equation}
\Re\partial_{0}=\frac{1}{2}\left(\partial_{0}+\partial_{0}^{*}\right)=\rho.\label{eq:d0posdef}
\end{equation}
This observation, or direct computation, gives that for bounded selfadjoint
positive $M_{0}:H\to H$, bounded $M_{1}:H\to H$ and skew-selfadjoint
$A:D\left(A\right)\subseteq H\to H$ 
\[
\Re\left\langle u|\left(\partial_{0}M_{0}+M_{1}+A\right)u\right\rangle _{\rho,0,0}=\left\langle u|\left(\rho M_{0}+\Re M_{1}\right)u\right\rangle _{\rho,0,0}
\]
for $u\in D\left(\partial_{0}\right)\cap D$$\left(A\right)$. With
the assumption that 
\begin{equation}
\rho M_{0}+\Re M_{1}\geq c_{0}>0\label{eq:well-pos}
\end{equation}
for all sufficiently large $\rho\in\left(0,\infty\right)$, we get
that the closure $\overline{\partial_{0}M_{0}+M_{1}+A}$ and its adjoint
$\overline{\partial_{0}M_{0}+M_{1}^{*}-A}=\left(\partial_{0}M_{0}+M_{1}+A\right)^{*}$
both have continuous inverses bounded by $\frac{1}{c_{0}}$. In particular,
the null space of $\overline{\partial_{0}M_{0}+M_{1}+A}$ and $\left(\partial_{0}M_{0}+M_{1}+A\right)^{*}$
are both trivial. Thus, we have the following well-posedness result,
see e.g. \cite{PIC_2010:1889,PDE_DeGruyter}.

\begin{thm} \label{thm:well} Let $M_{k}:H\to H${, $k=0,1$, be
continuous linear operators, $M_{0}$ selfadjoint, such that }\eqref{eq:well-pos}
holds for some $c_{0}\in\left(0,\infty\right)$ and for all $\rho\in\left(\rho_{0},\infty\right)$
with $\rho_{0}\in\left(0,\infty\right)$ sufficiently large. Moreover
let $A:D\left(A\right)\subseteq H\to H$ be skew-selfadjoint. Then
\[
\overline{\partial_{0}M_{0}+M_{1}+A}u=f
\]
has for any $f\in H_{\rho,0}\left(\mathbb{R},H\right)$ a unique solution
$u\in H_{\rho,0}\left(\mathbb{R},H\right)$. Furthermore, $u$ depends
on $f$ continuously, i.e. 
\[
\overline{\partial_{0}M_{0}+M_{1}+A}^{-1}:H_{\rho,0}\left(\mathbb{R},H\right)\to H_{\rho,0}\left(\mathbb{R},H\right)
\]
is a continuous linear operator{ for }$\rho\in\left(\rho_{0},\infty\right)$.\end{thm}

As a refinement of \eqref{eq:d0posdef} we find from integration by
parts that for $u\in\interior C_{1}\left(\mathbb{R},H\right)$ (and
so for $u\in D\left(\overline{\partial_{0}M_{0}+M_{1}+A}\right)$)
we have 
\[
\Re\left\langle u|\chi_{_{(-\infty,a]}}\left(\partial_{0}M_{0}+M_{1}\right)u\right\rangle _{\rho,0,0}\geq c_{0}\left\langle \chi_{_{(-\infty,a]}}u|\chi_{_{(-\infty,a]}}u\right\rangle _{\rho,0,0},
\]
since $|u(a)|^{2}e^{-2a\rho}\ge0$. This yields that we have also
causality in the sense of the following theorem.

\begin{thm}\label{thm:(Causality)}(Causality) Under the assumptions
of Theorem \ref{thm:well} we have 
\[
\chi_{_{(-\infty,a]}}\overline{\partial_{0}M_{0}+M_{1}+A}^{-1}=\chi_{_{(-\infty,a]}}\overline{\partial_{0}M_{0}+M_{1}+A}^{-1}\chi_{_{(-\infty,a]}}
\]
for all sufficiently large $\rho\in(0,\infty)$ .\end{thm}

We plan to approach the topic of boundary damping for the Timoshenko
beam within this abstract framework, which simplifies matters in so
far as we need only ensure that the spatial operator $A$ is skew-selfadjoint
and that assumptions of the type \eqref{eq:well-pos} hold for $M_{0}$,
$M_{1}$.

%%%%%%%%%%%%%%%%%%%%%%%%%%%%%%%%%%%%%%%%%%%

\section{Boundary Damping\label{sec-BD}}

For implementing suitable boundary conditions we consider first a
simple example of boundary conditions for the Timoshenko system. Assuming
that the beam is described by the unit interval $({-1/2},{+1/2})$,
which can always be achieved by translation and re-scaling, following
\cite{Roux2015194} we consider the case of the set of boundary conditions:
\begin{equation}
\begin{array}{rl}
V_{1}\left(\:\cdot\:,-1/2+0\right) & =0,\\
V_{1}\left(\:\cdot\:,1/2-0\right)+c\eta\left(\:\cdot\:,1/2-0\right) & =0,\\
s\left(\:\cdot\:,-1/2+0\right) & =0,\\
s\left(\:\cdot\:,1/2-0\right) & =0.
\end{array}\label{eq:bcs}
\end{equation}
To implement these boundary conditions we consider, in the spirit
of abstract grad-div systems as discussed in \cite{Picard20164888}
, the modified system 
\begin{equation}
\left(\partial_{0}M\left(\partial_{0}^{-1}\right)+A\right)\left(\begin{array}{c}
V_{1}\\
\left(\begin{array}{c}
\eta\\
\tau_{+}
\end{array}\right)\\
s\\
V_{2}
\end{array}\right)=\left(\begin{array}{c}
0\\
\left(\begin{array}{c}
f_{1}\\
0
\end{array}\right)\\
f_{2}\\
0
\end{array}\right),\label{eq:Timoshenko-damped-1}
\end{equation}
with material law of the simple form 
\begin{eqnarray*}
M\left(\partial_{0}^{-1}\right) & = & M_{0}+\partial_{0}^{-1}M_{1}
\end{eqnarray*}
with 
\[
M_{0}\coloneqq\left(\begin{array}{cccc}
\kappa_{1} & (0,0) & 0 & 0\\
\left(\begin{array}{cc}
0\\
0
\end{array}\right) & \left(\begin{array}{cc}
\nu_{1} & 0\\
0 & 0
\end{array}\right) & \left(\begin{array}{cc}
0\\
0
\end{array}\right) & \left(\begin{array}{cc}
0\\
0
\end{array}\right)\\
0 & (0,0) & \nu_{2} & 0\\
0 & (0,0) & 0 & \kappa_{2}
\end{array}\right),
\]
\[
M_{1}\coloneqq\left(\begin{array}{cccc}
0 & (0,0) & 0 & 0\\
\left(\begin{array}{cc}
0\\
0
\end{array}\right) & \left(\begin{array}{cc}
0 & 0\\
0 & c
\end{array}\right) & \left(\begin{array}{cc}
0\\
0
\end{array}\right) & \left(\begin{array}{cc}
1\\
0
\end{array}\right)\\
0 & (0,0) & d & 0\\
0 & (-1,0) & 0 & 0
\end{array}\right),
\]
where, we allow $d$ to be a continuous linear operator and $\nu_{1},\kappa_{1},\nu_{2},\kappa_{2}$
continuous selfadjoint and strictly positive definite operators in
$L^{2}\left({-1/2},{1/2}\right)$. For $c:\mathbb{C}\to\mathbb{C}$
we assume that $c$ is just multiplication by a positive real number.
Finally we set $A$ to be the skew-selfadjoint operator 
\[
A=\left(\begin{array}{cccc}
0 & B^{*} & 0 & 0\\
-B & \left(\begin{array}{cc}
0 & 0\\
0 & 0
\end{array}\right) & \left(\begin{array}{cc}
0\\
0
\end{array}\right) & \left(\begin{array}{cc}
0\\
0
\end{array}\right)\\
0 & (0,0) & 0 & -\partial_{1}\\
0 & (0,0) & -\interior{\partial}_{1} & 0
\end{array}\right).
\]
The differentiation operator $\partial_{1}$ denotes the weak derivative
in $L^{2}\left({-1/2},{1/2}\right)$ and the notation $\interior{\partial}_{1}$
indicates the use of Dirichlet boundary conditions, which has adjoint
$-\partial_{1}$ making $\interior{\partial}_{1}$ and $\partial_{1}$
skew-adjoint to each other.

Here we take $B:D\left(B\right)\subseteq L^{2}\left(-1/2,1/2\right)\to L^{2}\left(-1/2,1/2\right)\oplus\mathbb{C}\equiv\left(\begin{array}{c}
L^{2}\left(-1/2,1/2\right)\\
\mathbb{C}
\end{array}\right)$, where $D\left(B\right)\coloneqq\left\{ \varphi\in H_{1}\left(-1/2,1/2\right)|\varphi\left(-1/2+0\right)=0\right\} $,
given by 
\begin{equation}
B=\left(\begin{array}{c}
\overset{\triangleright}{\partial}_{1}\\
\delta_{\left\{ 1/2-0\right\} }
\end{array}\right).\label{eq:B}
\end{equation}
Here $\overset{\triangleright}{\partial}_{1}$ denotes $\partial_{1}$
with the domain constraint that it only acts on functions in $D\left(B\right)$,
which in particular vanish at $-1/2$.

%In the above we have 
%\begin{eqnarray*}
%\delta_{\left\{ 1/2-0\right\} }:H_{1}\left({-1/2},{1/2}\right) & \to\mathbb{C}\\
%\varphi & \mapsto\varphi\left(1/2-0\right).
%\end{eqnarray*}
%We note that $A$ is skew-selfadjoint by construction. 

The 1st, 3rd and 4th boundary conditions follow from the definitions
of the domains of $\interior{\partial}_{1}$ and $\overset{\triangleright}{\partial}_{1}$.
The first row of this system is the equation 
\[
B^{*}\left(\begin{array}{c}
\eta\\
\tau_{+}
\end{array}\right)+\partial_{0}(\kappa_{1}V_{1})=0.
\]
For $u\in D\left(B\right)$ with respect to the inner product $\left\langle \:\cdot\:|\:\cdot\:\right\rangle $
in $L^{2}(-1/2,1/2)$ we have 
\[
\left\langle u\Big|B^{*}\left(\begin{array}{c}
\eta\\
\tau_{+}
\end{array}\right)\right\rangle +\left\langle u|\partial_{0}(\kappa_{1}V_{1})\right\rangle =0
\]
so 
\[
\left\langle Bu\Big|\left(\begin{array}{c}
\eta\\
\tau_{+}
\end{array}\right)\right\rangle +\left\langle u|\partial_{0}(\kappa_{1}V_{1})\right\rangle =0
\]
i.e. 
\[
\left\langle \left(\begin{array}{c}
\partial_{1}u\\
u(1/2-0)
\end{array}\right)\Big|\left(\begin{array}{c}
\eta\\
\tau_{+}
\end{array}\right)\right\rangle +\left\langle u|\partial_{0}(\kappa_{1}V_{1})\right\rangle =0.
\]
Restricting to $u\in\interior{H}_{1}\left({-1/2},{1/2}\right)$ the
above equation yields 
\[
-\partial_{1}\eta+\partial_{0}(\kappa_{1}V_{1})=0.
\]
Now for $u\in\interior{H}_{1}\left({-1/2},{1/2}\right)$ with $u(-1/2+0)=0$,
integrating by parts we have 
\[
\overline{\tau}_{+}u(1/2-0)+\overline{\eta}(\cdot,1/2-0)u(1/2-0)-\overline{\eta}(\cdot,-1/2+0)u(-1/2+0)=\int_{-1/2}^{1/2}\overline{(\partial_{1}\eta-\partial_{0}(\kappa_{1}V_{1}))}u
\]
which gives 
\[
{\tau}_{+}+{\eta}(\cdot,1/2-0)=0.
\]
The second term of the second row in our system gives 
\[
-V_{1}(\cdot,1/2-0)+c\tau_{+}=0
\]
which combined with the previous expression gives 
\[
V_{1}(\cdot,1/2-0)+c{\eta}(\cdot,1/2-0)=0.
\]

The fact that system \eqref{eq:Timoshenko-damped-1} falls into the
class of abstract operators considered in \cite{Pi2009-1} now yields
the following well-posedness theorem. \begin{thm} Let $c$ be a positive
real number, $\nu_{1},\,\kappa_{1},\,\kappa_{2}$ continuous, selfadjoint,
strictly positive definite linear operators in $L^{2}\left({-1/2},{1/2}\right)$.
For the operators $\nu_{1}$ and $d$ we require 
\[
\rho\nu_{2}+\Re d\geq c_{0}>0
\]
for all sufficiently large $\rho\in(0,\infty)\:$. Then for any right-hand
side $F=\left(0,\left(f_{1},0\right),f_{2},0\right)\in H_{\rho,0}\left(\mathbb{R},H\right)$
with 
\[
H=L^{2}\left({-1/2},{1/2}\right)\oplus\left(L^{2}\left({-1/2},{1/2}\right)\oplus\mathbb{C}\right)\oplus L^{2}\left({-1/2},{1/2}\right)\oplus L^{2}\left({-1/2},{1/2}\right)
\]
there is a unique solution $\left(V_{1},\left(\eta,\tau_{+}\right),s,V_{2}\right)\in H_{\rho,0}\left(\mathbb{R},H\right)$.
\end{thm}

\begin{rem} In particular, the boundary conditions \eqref{eq:bcs}
are satisfied in the sense of $H_{\rho,-1}\left(\mathbb{R},\mathbb{C}\right)$,
which is the space of distributional temporal derivatives of $H_{\rho,0}\left(\mathbb{R},\mathbb{C}\right)$.
Moreover, $F=\left(f_{0},\left(f_{1},g_{1}\right),f_{2},f_{3}\right)\in H_{\rho,0}\left(\mathbb{R},H\right)$
can be completely arbitrary so that automatically an inhomogeneous
boundary condition of the form 
\[
V_{1}\left(\:\cdot\:,1/2-0\right)+c\eta\left(\:\cdot\:,1/2-0\right)=-g_{1}
\]
can be incorporated. Furthermore, the solution theory extends to data
in the space $H_{\rho,-\infty}\left(\mathbb{R},\mathbb{C}\right)$
of finite order distributional temporal derivatives of $H_{\rho,0}\left(\mathbb{R},\mathbb{C}\right)$,
see \cite{Pi2009-1} for details.\end{rem}

%%%%%%%%%%%%%%%%%%%%%%%%%%%%%%%%%%%%

\section{Other Material Laws.}

Another special case of boundary damping from the literature, allowing
for a non-vanishing right-hand side in \cite{Zietsman2004199}, is,
with adapted names of parameters and variables, given by replacing
the second equation of (\ref{eq:bcs}) with 
\[
\partial_{0}\tilde{I}\eta\left(\:\cdot\:,1/2-0\right)+V_{1}\left(\:\cdot\:,1/2-0\right)+{c}\eta\left(\:\cdot\:,1/2-0\right)=-g_{2}.
\]
This amounts to replacing $M_{0}$ in the above with 
\[
M_{0}\coloneqq\left(\begin{array}{cccc}
\kappa_{1} & (0,0) & 0 & 0\\
\left(\begin{array}{cc}
0\\
0
\end{array}\right) & \left(\begin{array}{cc}
\nu_{1} & 0\\
0 & \tilde{I}
\end{array}\right) & \left(\begin{array}{cc}
0\\
0
\end{array}\right) & \left(\begin{array}{cc}
0\\
0
\end{array}\right)\\
0 & (0,0) & \nu_{2} & 0\\
0 & (0,0) & 0 & \kappa_{2}
\end{array}\right).
\]
The parameter $c$ along with the additional parameter $\tilde{I}$
are now allowed to be non-negative reals with not both zero. In regards
to the model equations used in \cite{Roux2015194}, which differ slightly
from the above system of partial differential equations for the Timoshenko
beam, it may be advisable to consider the more general variant 
\[
M_{1}\coloneqq\left(\begin{array}{cccc}
0 & (0,0) & 0 & 0\\
\left(\begin{array}{cc}
0\\
0
\end{array}\right) & \left(\begin{array}{cc}
0 & 0\\
0 & c
\end{array}\right) & \left(\begin{array}{cc}
0\\
0
\end{array}\right) & \left(\begin{array}{cc}
-\sigma_{0}^{*}\\
0
\end{array}\right)\\
0 & (0,0) & d & 0\\
0 & (\sigma_{0},0) & 0 & 0
\end{array}\right),
\]
with a non-vanishing number (or operator) $\sigma_{0}$. It may be
noteworthy that for $d=0=c$ we have a system in which the norm of
\[
\left(\begin{array}{cccc}
\sqrt{\kappa_{1}} & \qquad0 & \qquad0 & \qquad0\\
0 & \qquad\left(\begin{array}{cc}
\sqrt{\nu_{1}} & 0\\
0 & \sqrt{\tilde{I}}
\end{array}\right) & \qquad0 & \qquad0\\
0 & \qquad0 & \qquad\sqrt{\nu_{2}} & \qquad0\\
0 & \qquad0 & \qquad0 & \qquad\sqrt{\kappa_{2}}
\end{array}\right)\left(\begin{array}{c}
V_{1}\\
\left(\begin{array}{c}
\eta\\
\tau_{+}
\end{array}\right)\\
s\\
V_{2}
\end{array}\right)
\]
is conserved, assuming a pure initial value problem. We remark that
in some model equations in the literature the rotational displacement
enters with the opposite sign. This is, however, just a unitarily
congruent version of our system obtained via the unitary transformation
matrix 
\[
\left(\begin{array}{cccc}
1 & \qquad0 & \qquad0 & \qquad0\\
0 & \qquad\left(\begin{array}{cc}
-1 & 0\\
0 & 1
\end{array}\right) & \qquad0 & \qquad0\\
0 & \qquad0 & \qquad1 & \qquad0\\
0 & \qquad0 & \qquad0 & \qquad1
\end{array}\right).
\]
This transformation changes only $B$ and $M_{1}$ slightly to 
\[
B=\left(\begin{array}{c}
-\overset{\triangleright}{\partial}_{1}\\
\delta_{\left\{ 1/2-0\right\} }
\end{array}\right)
\]
and 
\[
M_{1}\coloneqq\left(\begin{array}{cccc}
0 & (0,0) & 0 & 0\\
\left(\begin{array}{cc}
0\\
0
\end{array}\right) & \left(\begin{array}{cc}
0 & 0\\
0 & c
\end{array}\right) & \left(\begin{array}{cc}
0\\
0
\end{array}\right) & \left(\begin{array}{cc}
\sigma_{0}^{*}\\
0
\end{array}\right)\\
0 & (0,0) & d & 0\\
0 & (-\sigma_{0},0) & 0 & 0
\end{array}\right),
\]

Indeed, completely general material law operators $M\left(\partial_{0}^{-1}\right)$
can be handled in the same way, where $M$ is a bounded operator-valued
function, analytic in a ball of positive radius $r$ around the point
$r$ on the real axis. The crucial assumption is that the numerical
range of $\partial_{0}M\left(\partial_{0}^{-1}\right)$ is in the
right half plane and uniformly bounded away from the imaginary axis
for all sufficiently large $\rho\in(0,\infty)$:
\begin{equation}
\Re\left\langle u|\partial_{0}M\left(\partial_{0}^{-1}\right)u\right\rangle _{\rho,0,0}\geq c_{0}\left\langle u|u\right\rangle _{\rho,0,0}\label{eq:pos-def-general}
\end{equation}
for some $c_{0}\in\left(0,\infty\right)$ and all $u\in D\left(\partial_{0}\right)$.
This allows for a number of more intricate coupling phenomena.

If we assume that 
\[
M\left(\partial_{0}^{-1}\right)=\left(\begin{array}{cccc}
M_{00}\left(\partial_{0}^{-1}\right) & M_{01}\left(\partial_{0}^{-1}\right) & 0 & 0\\
M_{10}\left(\partial_{0}^{-1}\right) & M_{11}\left(\partial_{0}^{-1}\right) & 0 & 0\\
0 & 0 & M_{22}\left(\partial_{0}^{-1}\right) & M_{23}\left(\partial_{0}^{-1}\right)\\
0 & 0 & M_{32}\left(\partial_{0}^{-1}\right) & M_{33}\left(\partial_{0}^{-1}\right)
\end{array}\right)
\]
the system actually decomposes%
\footnote{Indeed, the coupling of all four equations is due to the off-diagonal
entries in $\mathcal{M}_{1}$ of the original first order Timoshenko
beam model.%
} into two $2\times2$-systems. Focusing on the first block system
we get 
\[
\left(\partial_{0}\left(\begin{array}{cc}
M_{00}\left(\partial_{0}^{-1}\right) & M_{01}\left(\partial_{0}^{-1}\right)\\
M_{10}\left(\partial_{0}^{-1}\right) & M_{11}\left(\partial_{0}^{-1}\right)
\end{array}\right)+\left(\begin{array}{cc}
0 & B^{*}\\
-B & 0
\end{array}\right)\right)\left(\begin{array}{c}
V_{1}\\
\left(\begin{array}{c}
\eta\\
\tau_{+}
\end{array}\right)
\end{array}\right)=\left(\begin{array}{c}
0\\
\left(\begin{array}{c}
f_{1}\\
0
\end{array}\right)
\end{array}\right)
\]
where $B$ is as in (\ref{eq:B}).

If we take 
\[
\tilde{M}(\partial_{0}^{-1})=\left(\begin{array}{cc}
M_{00}\left(\partial_{0}^{-1}\right) & M_{01}\left(\partial_{0}^{-1}\right)\\
M_{10}\left(\partial_{0}^{-1}\right) & M_{11}\left(\partial_{0}^{-1}\right)
\end{array}\right)=\left(\begin{array}{cc}
r+\partial_{0}^{-1}q & 0\\
0 & \left(\begin{array}{ccc}
s_{0}+\partial_{0}^{-1}s_{1} & 0 & 0\\
0 & \mu_{-}\left(\partial_{0}^{-1}\right) & 0\\
0 & 0 & \mu_{+}\left(\partial_{0}^{-1}\right)
\end{array}\right)
\end{array}\right)
\]
then for $s_{0}=p^{-1}$ and $s_{1}=0$ we obtain the undamped hyperbolic
case, while for $s_{0}=0$ and $s_{1}=p^{-1}$ we have the parabolic
case of the first order formulation of time-dependent Sturm-Liouville
problems (in standard terms). Here we replace $B$ by 
\[
\tilde{B}=\left(\begin{array}{c}
\partial_{1}\\
\delta_{\{-1/2+0\}}\\
\delta_{\{1/2-0\}}
\end{array}\right).
\]
The resulting operator equation 
\[
\left(\partial_{0}\tilde{M}(\partial_{0}^{-1})+\left(\begin{array}{cc}
0 & B^{*}\\
-B & 0
\end{array}\right)\right)\left(\begin{array}{c}
V_{1}\\
\left(\begin{array}{c}
\eta\\
\tau_{+}\\
\tau_{-}
\end{array}\right)
\end{array}\right)=\left(\begin{array}{c}
0\\
\left(\begin{array}{c}
f_{1}\\
0\\
0
\end{array}\right)
\end{array}\right)
\]
results in the first two boundary conditions of (\ref{eq:bcs}) being
replaced by dynamic boundary conditions %\begin{align*}
%\partial_{0}\mu_{+}\left(\partial_{0}^{^{-1}}\right)\tau_{+}-V_{1}\left(\:\cdot\:,1/2-0\right) & =h_{+}\\
%\partial_{0}\mu_{+}\left(\partial_{0}^{^{-1}}\right)\eta\left(\:\cdot\:,1/2-0\right)+V_{1}\left(\:\cdot\:,1/2-0\right) & =-h_{+}
%\end{align*}
\begin{align*}
\partial_{0}\mu_{+}\left(\partial_{0}^{^{-1}}\right)\eta\left(\:\cdot\:,1/2-0\right)+V_{1}\left(\:\cdot\:,1/2-0\right) & =0\\
\partial_{0}\mu_{-}\left(\partial_{0}^{^{-1}}\right)\eta\left(\:\cdot\:,-1/2+0\right)-V_{1}\left(\:\cdot\:,-1/2+0\right) & =0
\end{align*}
where $\mu_{+},\:\mu_{-}$ satisfy the assumption stated for $M$
in \eqref{eq:pos-def-general}. The point-wise evaluated Fourier-Laplace
transformed system (i.e. replacing $\partial_{0}$ by $z={\it i}\lambda+\rho$)
is discussed in \cite{BBW,Russ} for the parabolic case, i.e. $s_{0}=0$,
where $R_{1}\left(z\right)\coloneqq z\mu_{+}\left(\frac{1}{z}\right),\: R_{0}\left(z\right)\coloneqq z\mu_{-}\left(\frac{1}{z}\right)$
define rational Nevanlinna functions $R_{0}$ and $R_{1}$ .

\section{Other boundary conditions.}

A variety of other boundary conditions are also accessible via the
above approach. For example, implementing dynamic boundary conditions
in all unknowns and for both boundary parts we replace $A$ with 
\[
\tilde{A}=\left(\begin{array}{cccc}
0 & \tilde{B}^{*} & 0 & 0\\
-\tilde{B} & 0 & 0 & 0\\
0 & 0 & 0 & \tilde{B}^{*}\\
0 & 0 & -\tilde{B} & 0
\end{array}\right).
\]
The system now becomes 
\[
\left(\partial_{0}M\left(\partial_{0}^{-1}\right)+\tilde{A}\right)\left(\begin{array}{c}
V_{1}\\
\left(\begin{array}{c}
\eta\\
\tau_{0,-}\\
\tau_{0,+}
\end{array}\right)\\
s\\
\left(\begin{array}{c}
V_{2}\\
\tau_{1,-}\\
\tau_{1,+}
\end{array}\right)
\end{array}\right)=\left(\begin{array}{c}
0\\
\left(\begin{array}{c}
f_{1}\\
0\\
0
\end{array}\right)\\
0\\
\left(\begin{array}{c}
f_{2}\\
0\\
0
\end{array}\right)
\end{array}\right)
\]
If here again we assume that 
\[
M\left(\partial_{0}^{-1}\right)=\left(\begin{array}{cccc}
M_{00}\left(\partial_{0}^{-1}\right) & M_{01}\left(\partial_{0}^{-1}\right) & 0 & 0\\
M_{10}\left(\partial_{0}^{-1}\right) & M_{11}\left(\partial_{0}^{-1}\right) & 0 & 0\\
0 & 0 & M_{22}\left(\partial_{0}^{-1}\right) & M_{23}\left(\partial_{0}^{-1}\right)\\
0 & 0 & M_{32}\left(\partial_{0}^{-1}\right) & M_{33}\left(\partial_{0}^{-1}\right)
\end{array}\right)
\]
also this system decouples into two $2\times2$-systems. Focusing
on the first block system we get 
\begin{align*}
\left(\partial_{0}\tilde{M}(\partial_{0}^{-1})+\left(\begin{array}{cc}
0 & \tilde{B}^{*}\\
-\tilde{B} & 0
\end{array}\right)\right)\left(\begin{array}{c}
V_{1}\\
\left(\begin{array}{c}
\eta\\
\tau_{-}\\
\tau_{+}
\end{array}\right)
\end{array}\right)=\left(\begin{array}{c}
0\\
\left(\begin{array}{c}
f\\
0\\
0
\end{array}\right)
\end{array}\right).
\end{align*}
which results in dynamic boundary conditions as noted in the previous
section, likewise with the Sturm-Liouville case being a particular
application.

%we get -- compared to the discussion of the previous section of this
%decoupling case -- an additional boundary equation

%\begin{align*}
%\partial_{0}\mu_{-}\left(\partial_{0}^{^{-1}}\right)\tau_{-}-V_{1}\left(\:\cdot\:,-1/2+0\right) & =h_{-}\\
%\partial_{0}\mu_{-}\left(\partial_{0}^{^{-1}}\right)\eta\left(\:\cdot\:,-1/2+0\right)-V_{1}\left(\:\cdot\:,-1/2+0\right) & =h_{-}
%\end{align*}
%\[
%\partial_{0}\mu_{-}\left(\partial_{0}^{^{-1}}\right)=\frac{P_{-}\left(\partial_{0}^{-1}\right)}{-Q_{-}\left(\partial_{0}^{-1}\right)}=\frac{\partial_{0}^{r_{-}}P_{-}\left(\partial_{0}^{-1}\right)}{-\partial_{0}^{r_{-}}Q_{-}\left(\partial_{0}^{-1}\right)}=\frac{P_{0}\left(\partial_{0}\right)}{-Q_{0}\left(\partial_{0}\right)}=R_{0}\left(\partial_{0}\right).
%\]
%%%%% bibliography %%%%%%%%%%%%%%%%%%%%%%%%%%%%%


\begin{thebibliography}{10}
\bibitem{BBW}{P. A. Binding, P. J. Browne, B. A. Watson,} {{Sturm-Liouville
problems with boundary conditions rationally dependent on the eigenparameter,
II}, {\em J. Comp. Applied Math.,} 148:147--168, 2002.}

\bibitem{Pi2009-1} R.~Picard. \newblock {A Structural Observation
for Linear Material Laws in Classical Mathematical Physics.} \newblock
{\em {Math. Methods Appl. Sci.}}, 32(14):1768--1803, 2009.

\bibitem{Picard201354} R.~Picard. \newblock {Mother Operators
and their Descendants}. \newblock {\em Journal of Mathematical
Analysis and Applications}, 403(1):54--62, 2013.

\bibitem{PIC_2010:1889} R.~Picard. \newblock {An Elementary Hilbert
Space Approach to Evolutionary Partial Differential Equations}. \newblock
{\em Rend. Istit. Mat. Univ. Trieste}, 42 suppl.:185--204, 2010.

\bibitem{PDE_DeGruyter} R.~Picard and D.~F. McGhee. \newblock
{\em Partial Differential Equations: A unified Hilbert Space Approach},
volume~55 of {\em {De Gruyter Expositions in Mathematics}}.
\newblock {De Gruyter. Berlin, New York. 518 p.}, 2011.

\bibitem{Picard20164888} Rainer Picard, Stefan Seidler, Sascha Trostorff,
and Marcus Waurick. \newblock On abstract grad-div systems. \newblock
{\em Journal of Differential Equations}, 260(6):4888 -- 4917, 2016.

\bibitem{PIC-WAT} R.~Picard and B.A.~Watson \newblock {Evoluationary
problems involving Sturm-Liouville operators}. \newblock {\em
Operator Theory: Advance and Applications}, 236:393--406, 2013.

\bibitem{Roux2015194} A.~Roux, A.J. van~der Merwe, and N.F.J. van
Rensburg. \newblock {Elastic waves in a Timoshenko beam with boundary
damping.} \newblock {\em Wave Motion}, 57:194 -- 206, 2015.

\bibitem{Russ} {E. M. Russakovskii,} {{The matrix Sturm-Liouville
problem with spectral parameter in the boundary conditions. Algebraic
and operator aspects}, {\em Trans. Moscow Math. Soc.}, \textbf{57}
(1996), 159-184.}

\bibitem{Zietsman2004199} L.~Zietsman, N.F.J. van Rensburg, and
A.J. van~der Merwe. \newblock {A Timoshenko beam with tip body
and boundary damping.} \newblock {\em Wave Motion}, 39(3):199
-- 211, 2004.\end{thebibliography}
\end{document}